# Convergence of excursion point processes and its applications to functional limit theorems of Markov processes on a half-line

KOUJI YANO

*Department of Mathematics, Graduate School of Science, Kobe University, Kobe, Japan.
E-mail: kyano@math.kobe-u.ac.jp*

Invariance principles are obtained for a Markov process on a half-line with continuous paths on the interior. The domains of attraction of the two different types of self-similar processes are investigated. Our approach is to establish convergence of excursion point processes, which is based on Itô's excursion theory and a recent result on convergence of excursion measures by Fitzsimmons and the present author.

## 1. Introduction

A strong Markov process on $[0, \infty)$ with continuous paths on the interior $(0, \infty)$ is characterized as follows. Its generator $\mathcal{L}$ is an extension of $\mathcal{L}_m = \frac{\mathrm{d}}{\mathrm{d}m}\frac{\mathrm{d}}{\mathrm{d}x}$ on $(0, \infty)$ for a *speed measure* $\mathrm{d}m$ on $(0, \infty)$ under Feller's boundary condition ([5]), which is given by

$$r\mathcal{L}u(0) = \int_{(0,\infty)} \{u(x) - u(0)\}j(\mathrm{d}x) + cu'(0) \tag{1.1}$$

for constants $c, r \geq 0$ and a *jumping-in measure* $j$ on $(0, \infty)$. Itô and McKean [9] and Itô [8] have constructed a sample path of the strong Markov process characterized as above for a possible triplet $(m, j, c, r)$. Such a process which starts from the origin will be denoted by $X_{m,j,c,r}$.

Lamperti [18] has characterized the totality of strong Markov processes $X = X_{m,j,c,r}$ with the *self-similar property*, $(\lambda^{-\alpha}X(\lambda t): t \geq 0) \stackrel{\text{law}}{=} (X(t): t \geq 0)$ for some $\alpha > 0$. Such a process behaves as a Bessel diffusion on the interior $(0, \infty)$ and its behavior when it starts from the origin has the following two possibilities:

(a) *it enters the interior continuously, that is, it is a reflecting Bessel process;*

(b) *it jumps into the interior according to the jumping-in measure* $j = j^{(\beta)}$, *where* $j^{(\beta)}$ *is given by* $j^{(\beta)}((x,\infty)) = x^{-\beta}$.

The purpose of the present paper is to establish invariance principles for the process $X_{m,j,c,r}$. The domain of attraction for the possible limit process (a) or (b) varies according to whether the integral $\int^\infty xj(\mathrm{d}x)$ converges or diverges. The result is a generalization in our class of strong Markov processes of that of Stone [23], who has characterized the domain of attraction of the case (a) in the class of diffusion processes (without jumps at the origin). For this purpose, we appeal to the method of *convergence of excursion point processes* explained below, which enables us to understand clearly what happens in the excursion level. In the proof of our results, a crucial role is played by one of the main results of Fitzsimmons and Yano [6], who dealt with convergence of excursion measures for diffusion processes on $(0,\infty)$ via time-change of the Brownian excursion.

Let us give an example to illustrate the main theorems. Consider

$$\mathcal{L}_m = \frac{x}{2x+1} \frac{\mathrm{d}^2}{\mathrm{d}x^2} \qquad \text{on } (0,\infty), \tag{1.2}$$

that is, $\mathrm{d}m(x) = \frac{2x+1}{x}\mathrm{d}x$. The origin for $\mathcal{L}_m$ is exit but non-entrance and hence the continuous entrance is not allowed, namely, the constant $c$ must be 0. In particular, a reflecting $\mathcal{L}_m$-diffusion process does not exist. In addition, the process $X_{m,j,0,r}$ exists if and only if $\int_{0+} x\log(1/x)j(\mathrm{d}x) + j([1,\infty)) < \infty$ and either $r > 0$ or $j((0,1)) = \infty$ holds. Then, by Theorems 2.5 and 2.6, we obtain the following:

(i) *if $X_{m,j,0,r}$ is non-trivial and $\int^\infty xj(\mathrm{d}x) < \infty$, then the process $\frac{1}{\sqrt{\lambda}}X_{m,j,0,r}(\lambda \cdot)$ converges in law to a reflecting Brownian motion;*

(ii) *if $j((x,\infty)) \sim x^{-\beta}L(x)$ as $x \to \infty$ for $\beta \in (0,1)$ and some slowly varying function $L$ at infinity (with $\int^\infty xj(\mathrm{d}x) = \infty$ holding true in this case), then the process $\frac{1}{\sqrt{\lambda}}X_{m,j,0,r}(\lambda \cdot)$ converges in law to the process $X_{2x,j^{(\beta)},0,0}(\cdot)$.*

The method of the time-change of Brownian motion is quite useful to functional limit theorems of diffusion processes. For example, see [17, 19, 23, 24]. Recently, Fitzsimmons and Yano [6] have obtained limit theorems where the method of *the time-change of the Brownian excursion* is fully exploited. In the present paper, based on Itô's excursion theory ([7, 8]) and the method of the time-change of the Brownian excursion, we construct sample paths of the processes $X_{m,j,c,r}$ simultaneously for all possible characteristics $(m,j,c,r)$ from a common excursion point process. Our limit theorems are then reduced to certain continuity lemmas of $X_{m,j,c,r}$ and its inverse local time process $\eta_{m,j,c,r}$ with respect to $(m,j,c,r)$.

The key to our limit theorems is *convergence of excursion point processes*, which is stated in Propositions 4.2 and 4.3. Vague or other convergences of Poisson point processes on finite-dimensional spaces have been studied by many authors; see, for example, [3, 4, 10, 11, 12, 14, 15, 16]. For our purposes, we need a certain stronger convergence of Poisson point processes on the space of excursions. Let us roughly explain the idea. The excursion point process $\boldsymbol{N}_{m,j,c}$ of the process $X_{m,j,c,r}$ is realized as the image measure of a certain time-changed path $e_{m,j,c}$ under the excursion point process $\hat{\boldsymbol{N}}$ of a Brownian motion (see

Lemma 3.3). The propositions then assert that if $(m_\lambda, j_\lambda, c_\lambda)$ converges to $(m, j, c)$ in a certain sense, then $e_{m_\lambda, j_\lambda, c_\lambda}$ converges to $e_{m,j,c}$ in a certain sense for all points in the support of the excursion point process $\hat{\boldsymbol{N}}$ almost surely. The convergence $e_{m_\lambda, j_\lambda, c_\lambda} \to e_{m,j,c}$ under $\hat{\boldsymbol{N}}$ implies convergence of excursion point processes $\boldsymbol{N}_{m_\lambda, j_\lambda, c_\lambda} \to \boldsymbol{N}_{m,j,c}$. This may be regarded as an analogue of Skorokhod representation, which asserts that weak convergence of probability measures can be realized as almost-sure convergence of random variables on a certain probability space. We point out that our convergence of excursion point processes in the above sense is stronger than the vague convergence of those.

The present paper is organized as follows. In Section 2, we state the main theorems. In Section 3, we follow Itô [8] to construct a sample path of the process from an excursion point process. Continuity lemmas of excursion point processes which play important roles in proving our main theorems are stated in Section 4 and proved in Section 5. Under certain extra assumptions, we prove almost-sure continuity lemmas for the inverse local time processes in Section 6 and for the strong Markov processes considered in Section 7. In Section 8, we remove the extra assumptions and obtain in-probability continuity lemmas. We then conclude by completing the proof of our invariance principles.

## 2. Main theorems

Let $m\colon (0,\infty) \to (-\infty, \infty)$ be a right-continuous and strictly increasing function. For such $m$, we denote $\mathcal{L}_m = \frac{\mathrm{d}}{\mathrm{d}m}\frac{\mathrm{d}}{\mathrm{d}x}$. We always assume that $\int_{0+} x\, \mathrm{d}m(x) < \infty$, that is, that the origin for $\mathcal{L}_m$ is an exit boundary. There then exists an absorbing $\mathcal{L}_m$-diffusion process starting from $x > 0$, whose law will be denoted by $\boldsymbol{Q}_m^x$. If $m(0+)$ is finite, that is, the origin for $\mathcal{L}_m$ is exit and entrance, we denote by $\boldsymbol{n}_m$ the excursion measure away from the origin for the reflecting $\mathcal{L}_m$-diffusion process. For a Radon measure $j$ on $(0,\infty)$ and for non-negative constants $c$ and $r$, we denote by $X_{m,j,c,r}$, if it exists, a strong Markov process starting from the origin whose generator is an extension of $\mathcal{L}_m$ on $(0,\infty)$ and which is subject to Feller's boundary condition (1.1). The following theorem is due to Feller [5] and Itô [8].

**Theorem 2.1.** *Let $j$ be a Radon measure on $(0,\infty)$ and let $c$ and $r$ be non-negative constants. Then the process $X_{m,j,c,r}$ exists if and only if the following conditions* (C) *and* (C+) *hold:*

(C) *the pair $(m,j)$ satisfies*

$$j((x_0, \infty)) + \int_{(0,x_0)} j(\mathrm{d}x) \int_0^x m((y, x_0))\, \mathrm{d}y < \infty \qquad (2.1)$$

*for some $x_0 > 0$, and*

$$c = 0 \qquad \text{in the case where } m(0+) = -\infty; \qquad (2.2)$$

(C+)　$r > 0$ *in the case where* $c = 0$ *and* $j((0, x_0)) < \infty$.

*If the process exists, then its excursion measure away from the origin is described as*

$$\boldsymbol{n}_{m,j,c}(\Gamma) = \int_{(0,\infty)} j(\mathrm{d}x) \boldsymbol{Q}_m^x(\Gamma) + c\boldsymbol{n}_m(\Gamma). \tag{2.3}$$

We will denote by $L_{m,j,c,r}(t)$ a version of the local time at the origin, chosen so that

$$\boldsymbol{P}\left[\int_0^\infty \mathrm{e}^{-t}\,\mathrm{d}L_{m,j,c,r}(t)\right] = \frac{1}{C_{m,j,c,r}}, \tag{2.4}$$

where

$$C_{m,j,c,r} = r + \int_{(0,\infty)} (1 - \mathrm{e}^{-t}) \boldsymbol{n}_{m,j,c}(\zeta(e) \in \mathrm{d}t), \tag{2.5}$$

$\zeta(e)$ being the lifetime of an excursion path $e$. We will denote the right-continuous inverse of $L_{m,j,c,r}$ by $\eta_{m,j,c,r}$.

**Remark 2.2.** Theorem 2.1 has been obtained by Feller [5] in the case where $c$ is general but $m(0+)$ is finite, and by Itô [8] in the case where $m$ is general but $c = 0$. We can prove Theorem 2.1 in full generality in the same way as Itô [8], so we omit the proof.

**Remark 2.3.** The condition (2.1) always implies that $\int_{0+} x j(\mathrm{d}x) < \infty$. The converse also holds if $m(0+)$ is finite.

**Example 2.4.** Let us give typical examples of $m$ and $j$. For $\alpha > 0$, we define

$$m^{(\alpha)}(x) = \begin{cases} (1-\alpha)^{-1} x^{1/\alpha - 1}, & \text{if } 0 < \alpha < 1, \\ \log x, & \text{if } \alpha = 1, \\ -(\alpha-1)^{-1} x^{1/\alpha - 1}, & \text{if } \alpha > 1. \end{cases} \tag{2.6}$$

For $\beta > 0$, we define a Radon measure $j^{(\beta)}$ on $(0, \infty)$ by

$$j^{(\beta)}(\mathrm{d}x) = \beta x^{-\beta - 1}\,\mathrm{d}x. \tag{2.7}$$

According to Lamperti [18], Theorem 5.2, the totality of self-similar processes in the class of our strong Markov processes $X = X_{m,j,c,r}$ consists of the following two classes:

(a) $X = X_{m^{(\alpha)},0,c,0}$ for some $0 < \alpha < 1$ and $c > 0$. The process $X$ is then a reflecting Bessel process of dimension $2 - 2\alpha \in (0, 2)$. The process $X$ has the $\alpha$-self-similar property

$$(\lambda^{-\alpha} X(\lambda t) : t \geq 0) \stackrel{\text{law}}{=} (X(t) : t \geq 0). \tag{2.8}$$

In addition, its inverse local time process $\eta = \eta_{m^{(\alpha)},0,c,0}$ is an $\alpha$-stable subordinator which has the $1/\alpha$-self-similar property.

(b) $X = X_{m^{(\alpha)}, j^{(\beta)}, 0, 0}$ for some $\alpha > 0$ and $\beta \in (0, 1/\alpha)$. The process $X$ also has the $\alpha$-self-similar property. In addition, its inverse local time process $\eta = \eta_{m^{(\alpha)}, j^{(\beta)}, 0, 0}$ is an $\alpha\beta$-stable subordinator which has the $1/(\alpha\beta)$-self-similar property.

We equip the set of cadlag paths with Skorokhod's $J_1$-topology, following Lindvall [20]; see also [11] and [12]. For cadlag paths $w_\lambda$ and $w$, we say that $w_\lambda \to w$ $(J_1)$ if there exists a family of homeomorphisms of $[0, \infty)$ denoted by $\{\Lambda_\lambda : \lambda > 0\}$ such that

$$\lim_{\lambda \to \infty} \sup_{t \in [0,T]} |\Lambda_\lambda(t) - t| = 0 \qquad \text{for all } T > 0 \tag{2.9}$$

and

$$\lim_{\lambda \to \infty} \sup_{t \in [0,T]} |w_\lambda(\Lambda_\lambda(t)) - w(t)| = 0 \qquad \text{for all } T > 0. \tag{2.10}$$

Note that compact uniform convergence always implies convergence $(J_1)$ and that the converse holds if the limit is a continuous path on $[0, \infty)$.

Generally speaking, invariance principles require one of the following conditions to hold:

(M1) $m(x) \sim (1-\alpha)^{-1} x^{1/\alpha - 1} K(x)$ as $x \to \infty$ for some $\alpha \in (0, 1)$;
(M2) $m(\lambda x) - m(\lambda) \sim (\log x) K(\lambda)$ as $\lambda \to \infty$ for all $x > 0$;
(M3) $m(\infty) < \infty$ and $m(\infty) - m(x) \sim (\alpha - 1)^{-1} x^{1/\alpha - 1} K(x)$ as $x \to \infty$ for some $\alpha \in (1, \infty)$.

Here, $K(x)$ denotes a slowly varying function at infinity. For the conditions (M2) and (M3), see, for example, [17] and also [6]. For a certain technical reason, we need the following assumption, stronger than (M1)–(M3):

(M) $dm(x) = m'(x) dx$ on $(x_0, \infty)$ for some $x_0 > 0$, where $m'(x)$ is a non-negative locally bounded measurable function such that

$$m'(x) \sim \alpha^{-1} x^{1/\alpha - 2} K(x) \qquad \text{as } x \to \infty \tag{2.11}$$

and $m$ satisfies an integrability condition $\int_{0+} x \log \log(1/x) \, dm(x) < \infty$.

We say that $X = X_{m,j,c,r}$ is trivial if $j = 0$ and $c = 0$, which is equivalent to saying that $X(t) \equiv 0$; in fact, the process $X_{m,j,c,r}$ starts from the origin and does not jump in $(0, \infty)$ nor enter $(0, \infty)$ continuously. We now state the main theorems of the present paper.

**Theorem 2.5 (The convergent case).** *Assume that the process $X_{m,j,c,r}$ exists and is not trivial and that the condition* (M) *holds for $\alpha \in (0, 1)$ and for some slowly varying function $K(x)$ at infinity. Assume, in addition, that the following holds:*

(J1) $\int^\infty x j(dx) < \infty$.

*Let $u(\lambda) = \lambda^{1/\alpha} K(\lambda)$. It then holds that*

$$\frac{1}{\lambda} X_{m,j,c,r}(u(\lambda) \cdot) \to Y^{(\alpha)}(\cdot) \qquad (J_1) \qquad \text{in law} \tag{2.12}$$

as $\lambda \to \infty$, where $Y^{(\alpha)}$ stands for a reflecting Bessel process of dimension $2 - 2\alpha$.

**Theorem 2.6 (The divergent case).** *Assume that the process $X_{m,j,c,r}$ exists and that the condition* (M) *holds for $0 < \alpha < \infty$ and for some slowly varying function $K(x)$ at infinity. Assume, in addition, that the following holds:*

(J2) $j((x, \infty)) \sim x^{-\beta} L(x)$ *as $x \to \infty$ for some $0 < \beta < \min\{1, 1/\alpha\}$ and for some slowly varying function $L(x)$ at infinity.*

*Let $u(\lambda) = \lambda^{1/\alpha} K(\lambda)$. It then holds that*

$$\frac{1}{\lambda} X_{m,j,c,r}(u(\lambda)\cdot) \to X_{m^{(\alpha)}, j^{(\beta)}, 0, 0}(\cdot) \qquad (J_1) \qquad \text{in law} \tag{2.13}$$

*as $\lambda \to \infty$.*

## 3. Construction of a sample path

Based on the method of Itô [8] for constructing a sample path of the process $X_{m,j,c,r}$ under Feller's boundary condition (1.1), we shall give a realization of the processes on a common probability space. For the general excursion theory, see [7, 21] and also [2].

Let $E$ denote the set of continuous paths $e : [0, \infty) \to [0, \infty)$ such that if $e(t_0) = 0$ for some $t_0 > 0$, then $e(t) = 0$ for all $t > t_0$. We call $\zeta = \zeta(e) = \inf\{t > 0 : e(t) = 0\}$ the *lifetime* of a path $e \in E$. Here, we follow the usual convention that $\inf \varnothing = \infty$. We equip $E$ with a compact uniform topology and denote by $\mathcal{B}(E)$ its Borel $\sigma$-field. For $e \in E$, we denote the first hitting time to $a \geq 0$ by $\tau_a = \tau_a(e) = \inf\{t \geq 0 : e(t) = a\}$. In particular, $\tau_0(e) = 0$ if $e(0) = 0$. The supremum value is denoted by $M = M(e) = \sup_{t \geq 0} e(t)$. Under our notation, we note that $\{\tau_a < \infty\} = \{M \geq a\}$ on $\{\zeta < \infty\}$.

We recall the Brownian excursion measure. Let $\boldsymbol{n}_{\mathrm{BE}}$ denote the excursion measure away from the origin of a reflecting Brownian motion. That is, $\boldsymbol{n}_{\mathrm{BE}}$ is a $\sigma$-finite measure on $E$ such that

$$\boldsymbol{n}_{\mathrm{BE}}(e(t+\cdot) \in \Gamma) = \int_{(0,\infty)} \boldsymbol{Q}^x_{\mathrm{BM}}(\Gamma) \boldsymbol{P}^0_{\mathrm{3B}}(e(t) \in \mathrm{d}x) \qquad \text{for } t > 0 \text{ and } \Gamma \in \mathcal{B}(E), \tag{3.1}$$

where $\boldsymbol{Q}^x_{\mathrm{BM}}$ stands for the law on $E$ of an absorbing Brownian motion starting from $x > 0$ and $\boldsymbol{P}^0_{\mathrm{3B}}$ for that of a 3-dimensional Bessel process starting from 0 with the generator $\frac{1}{2}\frac{\mathrm{d}^2}{\mathrm{d}x^2} + \frac{1}{x}\frac{\mathrm{d}}{\mathrm{d}x}$. It is obvious that $\boldsymbol{n}_{\mathrm{BE}}(E \setminus E_1) = 0$, where

$$E_1 = \{e \in E : e(0) = 0, 0 < \zeta(e) < \infty\}. \tag{3.2}$$

Just as an almost everywhere Brownian path does, an almost everywhere excursion path with respect to the Brownian excursion measure has local times, which is precisely stated as follows.

**Theorem 3.1 (See, e.g., [2]).** *There exist a measurable functional $\ell\colon [0,\infty)\times [0,\infty)\times E \to [0,\infty)$ and a set $E_2 \in \mathcal{B}(E)$ with $\boldsymbol{n}_{\mathrm{BE}}(E\setminus E_2)=0$ such that, for every fixed $e\in E_2$, the function $\ell(t,x) = \ell(t,x,e)$ satisfies the following:*

  (i) *the function $[0,\infty)\times [0,\infty) \ni (t,x) \mapsto \ell(t,x)$ is jointly continuous;*
 (ii) *for any $x>0$, the function $t \mapsto \ell(t,x)$ is non-decreasing;*
(iii) *$\int_0^t 1_A(e(s))\,\mathrm{d}s = 2\int_A \ell(t,x)\,\mathrm{d}x$ holds for all $t\geq 0$ and $A\in\mathcal{B}([0,\infty))$.*

We remark that it follows from the occupation formula (iii) and the bi-continuity (i) that

$$\ell(t,x) = \lim_{\varepsilon\to 0+} \frac{1}{2\varepsilon} \int_0^t 1_{[x,x+\varepsilon)}(e(s))\,\mathrm{d}s \qquad \text{for } e\in E_2. \tag{3.3}$$

Moreover, we remark that $\ell(t,0)=0$ holds for $\boldsymbol{n}_{\mathrm{BE}}$-almost everywhere excursion path, whereas $\ell(t,0)>0$ for almost everywhere Brownian path.

Following [6], we introduce the time-change of the Brownian excursion. For a right-continuous strictly increasing function $m\colon (0,\infty) \to (-\infty,\infty)$ such that $\int_{0+} x\,\mathrm{d}m(x)<\infty$, we define a clock $A_m(t) = A_m[e](t)$ by $A_m(t) = \int_{(0,\infty)} \ell(t,x)\,\mathrm{d}m(x)$. Lemma 2.4 of [6], which we may call a version of Jeulin's lemma (see also [13] and [22]), says that $A_m(t)<\infty$ for $\boldsymbol{n}_{\mathrm{BE}}$-almost every excursion path. We now define a time-changed excursion path by $e_m(t) = e(A_m^{-1}(t))$ for $t\geq 0$. For $x>0$, we define a shifted path $\theta_x(e)$ by

$$\theta_x(e)(\cdot) = \begin{cases} e(\tau_x(e)+\cdot), & \text{if } M(e)>x, \\ 0(\cdot), & \text{if } M(e)\leq x, \end{cases} \tag{3.4}$$

where $0\in E$ is defined by $0(t)\equiv 0$. We define $e_{m,x}$ by the time-changed excursion path of $\theta_x(e)$, which coincides with the shifted path of $e_m$, namely,

$$e_{m,x} = (\theta_x(e))_m = \theta_x(e_m). \tag{3.5}$$

Then, fundamental to our method are the following identities (see the equalities (2.13) and (2.17) and Theorem 2.5 of [6]). For any $\Gamma\in\mathcal{B}(E)$ such that $0\notin\Gamma$,

$$\boldsymbol{Q}_m^x(\Gamma) = x\boldsymbol{n}_{\mathrm{BE}}(e_{m,x}\in\Gamma) \tag{3.6}$$

and

$$\boldsymbol{n}_m(\Gamma) = \boldsymbol{n}_{\mathrm{BE}}(e_m\in\Gamma). \tag{3.7}$$

*Remark 3.2.* If $m(0+)$ is finite, then the measure $\boldsymbol{n}_m$ is the excursion measure of the reflecting $\mathcal{L}_m$-diffusion process in the usual sense. Otherwise, $\boldsymbol{n}_m$ is *never* an excursion measure for any strong Markov process since $\int_{0+} t\boldsymbol{n}_m(\zeta(e)\in\mathrm{d}t)=\infty$. Nevertheless, the measure $\boldsymbol{n}_m$, which we call the *generalized excursion measure*, gives a useful tool to consider limit theorems involving the case where the origin for $\mathcal{L}_m$ is exit but non-entrance. See [25] and [6] for details.

Let $j$ be a Radon measure $j$ on $(0,\infty)$ such that $\int_{0+} xj(\mathrm{d}x) < \infty$ and let $c \geq 0$ be a constant. For a such pair $(j,c)$, we define a function $J(z)$ on $(0,\infty)$ by

$$J(z) = \inf\left\{x > 0 : c + \int_{(0,x]} yj(\mathrm{d}y) > z\right\}. \tag{3.8}$$

Let $J : (0,\infty) \to [0,\infty]$ be a right-continuous non-decreasing function such that $J(\infty) = \infty$. Conversely, if such a function $J$ is given, then we recover a pair $(j,c)$ by setting

$$j(\mathrm{d}x) = \frac{\mathrm{d}J^{-1}(x)}{x} \quad \text{and} \quad c = c(J) = \inf\{z > 0 : J(z) > 0\}, \tag{3.9}$$

where $J^{-1}$ is the right-continuous inverse of $J$: $J^{-1}(x) = \inf\{z > 0 : J(z) > x\}$. We always identify $(j,c)$ with $J$ in this way. Set $d(J) = \sup\{z > 0 : J(z) < \infty\}$. Then, $d(J) = c(J) + \int_{(0,\infty)} yj(\mathrm{d}y)$.

Based on the identities (3.6) and (3.7), we obtain the following.

**Lemma 3.3.** *Let $m : (0,\infty) \to (-\infty,\infty)$ be a right-continuous strictly increasing function and $J : (0,\infty) \to [0,\infty]$ be a right-continuous non-decreasing function such that $J(\infty) = \infty$. Then, for any non-negative measurable functional $F$ on $E$ such that $F(0) = 0$, the identity*

$$\int_{(0,d(J))\times E} F(e_{m,J(z)})\,\mathrm{d}z \otimes \boldsymbol{n}_{\mathrm{BE}}(\mathrm{d}e) = \int_E F(e)\boldsymbol{n}_{m,j,c}(\mathrm{d}e) \tag{3.10}$$

*holds, where $j$ and $c$ are given by (3.9) and*

$$\boldsymbol{n}_{m,j,c}(\Gamma) = \int_{(0,\infty)} j(\mathrm{d}x)\boldsymbol{Q}_m^x(\Gamma) + c\boldsymbol{n}_m(\Gamma). \tag{3.11}$$

**Proof.** We divide the domain of the integral into the two disjoint intervals as $(0,d(J)) = (0,c] \cup (c,d(J))$ and in the integral on $(c,d(J))$, we change the variables by $x = J(z)$. The left-hand side of (3.10) then becomes

$$\int_{(0,\infty)\times E} F(e_{m,x})xj(\mathrm{d}x) \otimes \boldsymbol{n}_{\mathrm{BE}}(\mathrm{d}e) + c\int_E F(e_m)\boldsymbol{n}_{\mathrm{BE}}(\mathrm{d}e). \tag{3.12}$$

Using the identities (3.6) and (3.7), we rewrite the above expression as

$$\int_{(0,\infty)\times E} F(e)j(\mathrm{d}x) \otimes \boldsymbol{Q}_m^x(\mathrm{d}e) + c\int_E F(e)\boldsymbol{n}_m(\mathrm{d}e), \tag{3.13}$$

which is exactly the right-hand side of (3.10). □

Let $\hat{\boldsymbol{N}}$ be a Poisson point process on $(0,\infty) \times (0,\infty) \times E$ with its characteristic measure $\mathrm{d}s \otimes \mathrm{d}z \otimes \boldsymbol{n}_{\mathrm{BE}}(\mathrm{d}e)$ defined on a probability space $(\hat{\Omega}, \hat{\mathcal{F}}, \hat{\boldsymbol{P}})$. Lemma 3.3 then asserts that the excursion point process corresponding to the excursion measure $\boldsymbol{n}_{m,j,c}$ can be realized by the law of $e_{m,J(z)}$ under $\hat{\boldsymbol{N}}(\mathrm{d}s \times \mathrm{d}z \times \mathrm{d}e)$. We define a process $\hat{\eta}_{m,J,r} = (\hat{\eta}_{m,J,r}(s))$ as

$$\hat{\eta}_{m,J,r}(s) = rs + \int_{(0,d(J)) \times E} \zeta(e_{m,J(z)}) \hat{\boldsymbol{N}}((0,s] \times \mathrm{d}z \times \mathrm{d}e). \tag{3.14}$$

Here, we note that

$$\zeta(e_{m,J(z)}) = \begin{cases} A_m(\zeta) - A_m(\tau_{J(z)}), & \text{on } \{M(e) > J(z)\}, \\ 0, & \text{on } \{M(e) \leq J(z)\}. \end{cases} \tag{3.15}$$

Under the identifications (3.8) and (3.9) between $(j,c)$ and $J$, the conditions (C) and (C+) of Theorem 2.1 stated in terms of $(m,j,c)$ are translated into those in terms of $(m,J)$ as follows:

(C) the pair $(m,J)$ satisfies

$$\int_{(z_0,d(J))} \frac{\mathrm{d}z}{J(z)} + \int_{(c(J),z_0)} \frac{\mathrm{d}z}{J(z)} \int_{c(J)}^{J(z)} m((y,J(z_0))) \, \mathrm{d}y < \infty \tag{3.16}$$

for some $z_0 \in (c(J), d(J))$ and

$$c(J) = 0 \quad \text{in the case where } m(0+) = -\infty; \tag{3.17}$$

(C+) $r > 0$ in the case where $\int_{0+} \frac{\mathrm{d}z}{J(z)} < \infty$.

We then obtain

**Lemma 3.4.** *Let $m:(0,\infty) \to (-\infty,\infty)$ be a right-continuous strictly increasing function and $J:(0,\infty) \to [0,\infty]$ a right-continuous non-decreasing function such that $J(\infty) = \infty$. The process $\hat{\eta}_{m,J,r}$ is then a non-decreasing Lévy process if and only if the condition (C) holds. In this case, the Lévy measure is given by $\boldsymbol{n}_{m,j,c}(\zeta(e) \in \cdot)$. Moreover, the process $\hat{\eta}_{m,J,r}$ is increasing if and only if the condition (C+) holds.*

**Proof.** It is immediate by construction that $\hat{\eta}_{m,J,r} = (\hat{\eta}_{m,J,r}(s) : s \geq 0)$ is a Lévy process. The Laplace transform $\hat{\boldsymbol{P}}[e^{-\xi \hat{\eta}_{m,J,r}(s)}]$ is given by

$$\exp\left\{-\xi rs - s\int_{(0,\infty) \times E}(1 - e^{-\xi \zeta(e_{m,J(z)}) 1_{\{M(e) > J(z)\}}}) \, \mathrm{d}z \otimes \boldsymbol{n}_{\mathrm{BE}}(\mathrm{d}e)\right\}. \tag{3.18}$$

Using Lemma 3.3, we rewrite the expression (3.18) as

$$\exp\left\{-\xi rs - s\int_{(0,\infty)}(1 - e^{-\xi t}) \boldsymbol{n}_{m,j,c}(\zeta(e) \in \mathrm{d}t)\right\}. \tag{3.19}$$

It is well known that the integral

$$\int_{(0,\infty)} j(\mathrm{d}x) \int_{(0,\infty)} (1 - e^{-\xi t}) \boldsymbol{Q}_m^x(\zeta(e) \in \mathrm{d}t) \tag{3.20}$$

is finite for all $\xi > 0$ if and only if (2.1) (or (3.16)) of the condition (C) holds and that the integral

$$\int_{(0,\infty)} (1 - e^{-\xi t}) \boldsymbol{n}_m(\zeta(e) \in dt) \tag{3.21}$$

is finite for all $\xi > 0$ if and only if (2.2) (or (3.17)) of the condition (C) holds. Hence, we conclude that the condition (C) is the necessary and sufficient condition for the process $\hat{\eta}_{m,J,r}$ to be a Lévy process. It is obvious that the Lévy process $\hat{\eta}_{m,J,r}$ is strictly increasing if and only if the condition (C+) is satisfied. □

Suppose that the conditions (C) and (C+) hold. We define a process $\hat{X}_{m,J,r} = (\hat{X}_{m,J,r}(t))$ by setting

$$\hat{X}_{m,J,r}(t) = e_{m,J(z)}(t - \hat{\eta}_{m,J,r}(s-)) \tag{3.22}$$

if $\hat{\eta}_{m,J,r}(s-) \leq t < \hat{\eta}_{m,J,r}(s)$ for some point $(s,z,e)$ in the support of $\hat{\boldsymbol{N}}(ds \times dz \times de)$ and by setting $\hat{X}_{m,J,r}(t) = 0$ otherwise. We now have the following.

**Proposition 3.5.** *Let* $m:(0,\infty) \to (-\infty,\infty)$ *be a right-continuous strictly increasing function and* $J:(0,\infty) \to [0,\infty]$ *a right-continuous non-decreasing function such that* $J(\infty) = \infty$. *Suppose that conditions* (C) *and* (C+) *hold. The law of* $(\hat{X}_{m,J,r}, \hat{\eta}_{m,J,r})$ *on the probability space* $(\hat{\Omega}, \hat{\mathcal{F}}, \hat{\boldsymbol{P}})$ *is then identical to that of* $(X_{m,j,c,r}, \eta_{m,j,c,r})$.

The proof is obvious by Theorem 2.1 and Lemma 3.3, so we omit it. Therefore, we obtain a realization of the process $X_{m,j,c,r}$ defined on the common probability space $(\hat{\Omega}, \hat{\mathcal{F}}, \hat{\boldsymbol{P}})$.

**Remark 3.6.** If $m(0+)$ is finite, that is, the origin for $\mathcal{L}_m$ is exit and entrance, then the process $X_{m,0,c,0}$ for positive $c$ exists, which is exactly a reflecting $\mathcal{L}_m$-diffusion process starting from the origin. In this case, the function $J(z)$ is given by

$$J(z) = V_{(0,c)}(z) = \begin{cases} 0, & \text{for } 0 < z < c, \\ \infty, & \text{for } z \geq c. \end{cases} \tag{3.23}$$

## 4. Convergence of excursion point processes

For a function $m$ which satisfies either one of the three conditions (M1)–(M3), we set

$$m_\lambda(x) = \begin{cases} m(\lambda x)/\{\lambda^{1/\alpha - 1} K(\lambda)\}, & \text{if } 0 < \alpha < 1, \\ \{m(\lambda x) - m(\lambda)\}/\{\lambda^{1/\alpha - 1} K(\lambda)\}, & \text{if } \alpha = 1, \\ \{m(\lambda x) - m(\infty)\}/\{\lambda^{1/\alpha - 1} K(\lambda)\}, & \text{if } \alpha > 1 \end{cases} \tag{4.1}$$

so that $dm_\lambda(x) = dm(\lambda x)/\{\lambda^{1/\alpha - 1} K(\lambda)\}$ in all cases.

The following lemma plays an important role in the proofs of Theorems 2.5 and 2.6.

**Lemma 4.1.** *Let $v(\lambda)$ be an arbitrary function. The identity in law*

$$\left(\frac{1}{\lambda}\hat{X}_{m,J,r}(u(\lambda)\cdot), \frac{1}{u(\lambda)}\hat{\eta}_{m,J,r}(v(\lambda)\cdot)\right) \stackrel{\text{law}}{=} (\hat{X}_{m_\lambda,J_\lambda,r_\lambda}(\cdot), \hat{\eta}_{m_\lambda,J_\lambda,r_\lambda}(\cdot)) \tag{4.2}$$

*holds, where*

$$r_\lambda = \frac{v(\lambda)}{u(\lambda)}r \tag{4.3}$$

*and $J_\lambda$ is defined by*

$$J_\lambda(z) = \frac{1}{\lambda}J\left(\frac{\lambda}{v(\lambda)}z\right). \tag{4.4}$$

**Proof.** For $e \in E$ and $\lambda > 0$, we define $e^\lambda \in E$ by $e^\lambda(t) = \lambda e(t/\lambda^2)$. Then, $\boldsymbol{n}_{\mathrm{BE}}(e^\lambda \in \cdot) = \lambda \boldsymbol{n}_{\mathrm{BE}}(\cdot)$ and we hence obtain

$$\left\{\int 1_A(s,x,e)\hat{\boldsymbol{N}}(\mathrm{d}s \times \mathrm{d}x \times \mathrm{d}e) : A \in \mathcal{B}((0,\infty) \times (0,\infty) \times E)\right\}$$
$$\stackrel{\text{law}}{=} \left\{\int 1_A(v(\lambda)s, \lambda x/v(\lambda), e^\lambda)\hat{\boldsymbol{N}}(\mathrm{d}s \times \mathrm{d}x \times \mathrm{d}e) : A \in \mathcal{B}((0,\infty) \times (0,\infty) \times E)\right\}. \tag{4.5}$$

Using this identity in law, we immediately obtain (4.2). □

By the definition (4.1), it is immediate that

$$\lim_{\lambda \to \infty} m_\lambda(x) = m^{(\alpha)}(x) \qquad \text{for all } x > 0. \tag{4.6}$$

Consider the case of Theorem 2.5. Since $d(J) = c + \int_{(0,\infty)} yj(\mathrm{d}y)$, the assumption (J1) is equivalent to $d(J) < \infty$. We take $v(\lambda) = \lambda$ and adopt the notation of Lemma 4.1. We then see that

$$r_\lambda = \frac{r}{\lambda^{1/\alpha-1}K(\lambda)} \to 0 \quad \text{and} \quad J_\lambda(z) = J(z)/\lambda \to V_{(0,\tilde{c})}(z), \tag{4.7}$$

where $\tilde{c} = d(J) = c + \int_{(0,\infty)} yj(\mathrm{d}y)$. Here, the function $V_{(0,\tilde{c})}$ was introduced in (3.23).

Consider the case of Theorem 2.6. The assumption (J2) is equivalent to $J^{-1}(x) \sim \frac{\beta}{1-\beta}x^{1-\beta} \times L(x)$ as $x \to \infty$. We take $v(\lambda) = \lambda^\beta/L(\lambda)$ and adopt the notation of Lemma 4.1. We then see that

$$r_\lambda = \frac{r}{\lambda^{1/\alpha-\beta}K(\lambda)L(\lambda)} \to 0 \quad \text{and} \quad J_\lambda(z) = \frac{J(\lambda^{1-\beta}L(\lambda)z)}{\lambda} \to J^{(\beta)}(z), \tag{4.8}$$

where $J^{(\beta)}(z) = (\frac{1-\beta}{\beta}z)^{1/(1-\beta)}$, $c(J^{(\beta)}) = 0$ and $\mathrm{d}(J^{(\beta)})^{-1}(x) = xj^{(\beta)}(\mathrm{d}x)$.

Now, we may think that our problem is reduced to a suitable continuity of the excursion path $e_{m,J(z)}$ and of its lifetime $\zeta(e_{m,J(z)})$ with respect to $(m,J)$ for fixed points $(z, e)$. Central to our method are the following two continuity lemmas of excursion point processes.

**Proposition 4.2 (The convergent case).** *Suppose that* $\int_{0+} x \log\log(1/x)\,\mathrm{d}m(x) < \infty$ *and that any one of the three conditions* (M1), (M2) *and* (M3) *holds. Suppose, in addition, that the condition* (J1) *holds. Set* $v(\lambda) = \lambda$ *and adopt the notation (4.3) and (4.4). The following then holds with* $\hat{\boldsymbol{P}}$-*probability one:*

$$\lim_{\lambda \to \infty} \zeta(e_{m_\lambda, J_\lambda(z)}) = \zeta(e_{m^{(\alpha)}}). \tag{4.9}$$

*Further,*

$$\lim_{\lambda \to \infty} \sup_{t \geq 0} |e_{m_\lambda, J_\lambda(z)}(t) - e_{m^{(\alpha)}}(t)| = 0 \tag{4.10}$$

*holds for all* $(z, e)$ *in the support of the measure* $\hat{\boldsymbol{N}}((0, \infty) \times \mathrm{d}z \times \mathrm{d}e)$.

**Proposition 4.3 (The divergent case).** *Suppose that* $\int_{0+} x \log\log(1/x)\,\mathrm{d}m(x) < \infty$ *and that any one of the three conditions* (M1), (M2) *and* (M3) *holds. Suppose, in addition, that the condition* (J2) *holds. Set* $v(\lambda) = \lambda^\beta/L(\lambda)$ *and adopt the notation (4.3) and (4.4). The following then holds with* $\hat{\boldsymbol{P}}$-*probability one:*

$$\lim_{\lambda \to \infty} \zeta(e_{m_\lambda, J_\lambda(z)}) = \zeta(e_{m^{(\alpha)}, J^{(\beta)}(z)}). \tag{4.11}$$

*Further,*

$$\lim_{\lambda \to \infty} \sup_{t \geq 0} |e_{m_\lambda, J_\lambda(z)}(t) - e_{m^{(\alpha)}, J^{(\beta)}(z)}(t)| = 0 \tag{4.12}$$

*holds for all* $(z, e)$ *in the support of the measure* $\hat{\boldsymbol{N}}((0, \infty) \times \mathrm{d}z \times \mathrm{d}e)$.

The proofs of Propositions 4.2 and 4.3 will be given in the next section.

## 5. Proof of the continuity lemmas of the excursion point processes

We introduce the following assumption.

(A1) $m_\lambda(x) \to m_\infty(x)$ as $\lambda \to \infty$ for all continuity points $x > 0$ of $m_\infty$ and

$$\lim_{\delta \to 0+} \limsup_{\lambda \to \infty} \int_{(0,\delta]} x \log\log(1/x)\,\mathrm{d}m_\lambda(x) = 0. \tag{5.1}$$

Condition (5.1) is called $\mathcal{M}_L$-*tightness* in [6]. The following theorem plays a crucial role in the proof of our main theorems.

**Theorem 5.1 (Theorem 2.9 of [6]).** *Suppose that condition* (A1) *holds. Set*

$$E_3 = \left\{ e \in E : \lim_{\lambda \to \infty} \sup_{t \geq 0} |A_{m_\lambda}(t) - A_{m_\infty}(t)| = 0 \right\} \quad (5.2)$$

*and*

$$E_4 = \left\{ e \in E : \lim_{\lambda \to \infty} \sup_{t \geq 0} |e_{m_\lambda}(t) - e_{m_\infty}(t)| = 0 \right\}. \quad (5.3)$$

*Then,* $\boldsymbol{n}_{\mathrm{BE}}(E \setminus (E_3 \cap E_4)) = 0$.

For later use, we set $E^* = E_1 \cap E_2 \cap E_3 \cap E_4$ so that $\boldsymbol{n}_{\mathrm{BE}}(E \setminus E^*) = 0$.
In addition, we introduce the following assumption.

(A2) $J_\lambda(z) \to J_\infty(z)$ as $\lambda \to \infty$ for all $z > 0$ and the right-continuous inverse $J_\infty^{-1}(x) = \inf\{z > 0 : J_\infty(z) > x\}$ is absolutely continuous on $(0, \infty)$ with respect to the Lebesgue measure $\mathrm{d}x$.

Under these assumptions, we obtain the following.

**Lemma 5.2.** *Suppose that the conditions* (A1) *and* (A2) *hold. The following statement then holds with* $\hat{\boldsymbol{P}}$-*probability one:*

$$\lim_{\lambda \to \infty} \zeta(e_{m_\lambda, J_\lambda(z)}) = \zeta(e_{m_\infty, J_\infty(z)}). \quad (5.4)$$

*Further,*

$$\lim_{\lambda \to \infty} \sup_{t \geq 0} |e_{m_\lambda, J_\lambda(z)}(t) - e_{m_\infty, J_\infty(z)}(t)| = 0 \quad (5.5)$$

*holds for all* $(z, e)$ *in the support of* $\hat{\boldsymbol{N}}((0, \infty) \times \mathrm{d}z \times \mathrm{d}e)$.

**Proof.** Set

$$U = \left\{ (x, e) \in [0, \infty) \times E^* : \lim_{\varepsilon \to 0} \tau_{x+\varepsilon}(e) = \tau_x(e) \right\}. \quad (5.6)$$

Recall the definitions (5.2) and (5.3) and the identity (3.15). Then, by assumption (A2), it is obvious that the convergences (5.4) and (5.5) hold if $(J_\infty(z), e) \in U$. Hence, the desired convergence follows if we prove that, with $\hat{\boldsymbol{P}}$-probability one, the set

$$\{(z, e) \in (0, \infty) \times E : (J_\infty(z), e) \notin U\} \quad (5.7)$$

has null measure with respect to the point measure $\hat{\boldsymbol{N}}((0, \infty) \times \mathrm{d}z \times \mathrm{d}e)$. For this, it suffices to show that the set (5.7) has null measure with respect to the characteristic measure $\mathrm{d}z \otimes \boldsymbol{n}_{\mathrm{BE}}(\mathrm{d}e)$.

We note that

$$\lim_{x \to 0+} \tau_x(e) = 0 \quad \text{on } E_1 = \{e(0) = 0, 0 < \zeta(e) < \infty\}. \quad (5.8)$$

In fact, $\tau_x(e)$ converges decreasingly to some $t_0 \in [0, \zeta(e))$ as $x$ tends decreasingly to 0 and hence $x = e(\tau_x) \to e(t_0) = 0$ by the continuity of $e(t)$ at $t = 0$, which shows that $t_0 = 0$. Hence, we obtain $\boldsymbol{n}_{\mathrm{BE}}(\lim_{x \to 0+} \tau_x(e) \neq 0) = 0$, which shows that the set (5.7) restricted to $\{(z,e) : J_\infty(z) = 0\}$ has null measure with respect to the characteristic measure $\mathrm{d}z \otimes \boldsymbol{n}_{\mathrm{BE}}(\mathrm{d}e)$.

Let $e \in E$ be fixed for the time being. Since the function $(0, M(e)] \ni x \mapsto \tau_x(e)$ is non-decreasing and since $\tau_x(e) = \infty$ for all $x > M(e)$, we have $\lim_{\varepsilon \to 0} \tau_{x+\varepsilon}(e) = \tau_x(e)$ for $\mathrm{d}x$-almost every $x$. Hence, we conclude that the set (5.7) restricted to $\{(z,e) : 0 < J_\infty(z) < \infty\}$ has null measure with respect to the characteristic measure $\mathrm{d}z \otimes \boldsymbol{n}_{\mathrm{BE}}(\mathrm{d}e)$. The proof is now complete. $\square$

Let us reduce Propositions 4.2 and 4.3 to Lemma 5.2. For this, we check that assumptions (A1) and (A2) hold under each of the assumptions of Propositions 4.2 and 4.3.

The following lemma is a slight improvement of [6], Lemma 2.17.

**Lemma 5.3.** *If $\int_{0+} x \log \log(1/x) \, \mathrm{d}m(x) < \infty$ and if any one of the three conditions* (M1), (M2) *and* (M3) *is satisfied, then condition* (A1) *is satisfied.*

**Proof.** It is obvious that $m_\lambda(x) \to m(x)$ as $\lambda \to \infty$ for all $x > 0$. Hence, we need only to check condition (5.1). Set $m^0(x) = m(\max\{x, 1\})$ and $m^1(x) = m(\min\{x, 1\})$. Since $\mathrm{d}m(x) = \mathrm{d}m^0(x) + \mathrm{d}m^1(x)$, it suffices to show that the condition (5.1) is satisfied for $m = m^0$ and $m = m^1$.

For the proof of (5.1) for $m = m^0$, the same argument as used in [6], Lemma 2.17, is still valid and hence we omit it.

Writing $\mathrm{d}m^1_\lambda(x) = \mathrm{d}m(\lambda x)/\{\lambda^{1/\alpha - 1} K(\lambda)\}$, we have

$$\int_{(0,\delta]} x \log \log(1/x) \, \mathrm{d}m^1_\lambda(x) = \frac{1}{\lambda^{1/\alpha} K(\lambda)} \int_{(0,1]} 1_{(0,\lambda\delta]}(x) x \log \log(\lambda/x) \, \mathrm{d}m(x). \quad (5.9)$$

Using the inequality $a + b \leq (1+a)(1+b)$ for $a, b > 0$, we see that the right-hand side of (5.9) is dominated by

$$\frac{1 + \log\{1 + \log \lambda\}}{\lambda^{1/\alpha} K(\lambda)} \int_{(0,1]} x\{1 + \log\{1 + \log(1/x)\}\} \, \mathrm{d}m(x). \quad (5.10)$$

Since $\int_{0+} x \log \log(1/x) \, \mathrm{d}m(x) < \infty$, the integral in (5.10) converges and hence the expression (5.10) converges to zero as $\lambda \to \infty$, which shows that (5.1) holds for $m = m^1$. $\square$

*Remark 5.4.* Thanks to Lemma 5.3, some of the assumptions of [6], Theorem 2.16, can be relaxed—the assumption on $m$ near the origin can be replaced by the assumption $\int_{0+} x \times \log \log(1/x) \, \mathrm{d}m(x) < \infty$.

**Lemma 5.5.** *Suppose that either one of the conditions* (J1) *and* (J2) *holds. The condition* (A2) *is then satisfied, where $J_\infty = V_{(0,d(J))}$ in the former case and where $J_\infty = J^{(\beta)}$ in the latter case.*

**Proof.** This is immediate by (4.7) and (4.8). □

Combining Lemma 5.2 with Lemmas 5.3 and 5.5, we have completed the proofs of Propositions 4.2 and 4.3.

## 6. Convergence of the inverse local time process

The following two propositions, although they need extra assumptions, play an essential role in our proof of Theorems 2.5 and 2.6.

**Proposition 6.1 (The convergent case).** *Suppose that the conditions* (M) *and* (J1) *hold. Suppose, in addition, that* $dm(x)$ *has a locally bounded density on the whole of* $(0, \infty)$ *such that*

$$\limsup_{x \to 0+} \frac{m'(x)}{x^{1/\alpha - 2}} < \infty. \tag{6.1}$$

*Set* $v(\lambda) = \lambda$ *and adopt the notation (4.3) and (4.4). Then, with* $\hat{P}$*-probability one,*

$$\lim_{\lambda \to \infty} \sup_{s \in [0,S]} |\hat{\eta}_{m_\lambda, J_\lambda, r_\lambda}(s) - \hat{\eta}_{m^{(\alpha)}, V_{(0,d(J))}, 0}(s)| = 0 \qquad \text{for all } S > 0. \tag{6.2}$$

**Proposition 6.2 (The divergent case).** *Suppose that the conditions* (M) *and* (J2) *hold. Suppose, in addition, that* $dm(x)$ *has a locally bounded density on the whole of* $(0, \infty)$ *such that (6.1) holds and that*

$$\limsup_{x \to 0+} x^{\beta - 1} \int_{(0,x]} yj(dy) < \infty. \tag{6.3}$$

*Set* $v(\lambda) = \lambda^\beta / L(\lambda)$ *and adopt the notation (4.3) and (4.4). Then, with* $\hat{P}$*-probability one,*

$$\lim_{\lambda \to \infty} \sup_{s \in [0,S]} |\hat{\eta}_{m_\lambda, J_\lambda, r_\lambda}(s) - \hat{\eta}_{m^{(\alpha)}, J^{(\beta)}, 0}(s)| = 0 \qquad \text{for all } S > 0. \tag{6.4}$$

*Remark 6.3.* The conclusions (6.2) and (6.4) are uniform convergence instead of $J_1$-convergence, in spite of cadlag processes. The reason is that the processes involved jump at the same points.

Consider the following conditions:

(A3) $r_\lambda \to r_\infty$ as $\lambda \to \infty$;
(A4) each $m_\lambda$ has a locally bounded density, that is, $dm_\lambda(x) = m'_\lambda(x) dx$, and $m'_\lambda(x) \leq m'_+(x)$ and $J_\lambda(z) \geq J_+(z)$ hold for all $x, z \in (0, \infty)$ and $\lambda > 0$ for some $(m_+, J_+)$ which satisfies the conditions (3.16) and (3.17).

We then obtain the following continuity lemma of the Lévy process.

**Lemma 6.4.** *If the conditions* (A1)–(A4) *are satisfied, then the convergence*

$$\lim_{\lambda \to \infty} \sup_{s \in [0,S]} |\hat{\eta}_{m_\lambda, J_\lambda, r_\lambda}(s) - \hat{\eta}_{m_\infty, J_\infty, r_\infty}(s)| = 0 \tag{6.5}$$

*holds with $\hat{P}$-probability one for all $S > 0$.*

**Proof.** Recall that $\sup_{s \in [0,S]} |\hat{\eta}_{m_\lambda, J_\lambda, r_\lambda}(s) - \hat{\eta}_{m_\infty, J_\infty, r_\infty}(s)|$ is dominated by the sum of $|r_\lambda - r_\infty|$ and the integral $I_\lambda := \int_{(0,\infty) \times E} |F_\lambda(z,e) - F_\infty(z,e)| \hat{N}((0,S] \times dz \times de)$, where $F_\lambda(z,e) = \zeta(e_{m_\lambda, J_\lambda(z)}) 1_{\{M(e) > J_\lambda(z)\}}$ for $\lambda \leq \infty$. Set $F_+(z,e) = \zeta(e_{m_+, J_+(z)}) 1_{\{M(e) > J_+(z)\}}$. Since the variable $F_+(z,e)$ is integrable with respect to the measure $dz \otimes n_{\mathrm{BE}}(de)$, there exists $\hat{\Omega}^* \in \hat{\mathcal{F}}$ with $\hat{P}(\hat{\Omega}^*) = 1$ on which the variable $F_+(z,e)$ is integrable with respect to the measure $\hat{N}((0,S] \times dz \times de)$ and $\hat{N}((0,S] \times (0,\infty) \times (E \setminus E^*)) = 0$. Let $\hat{\omega}^* \in \hat{\Omega}^*$ be fixed.

By the conditions (A1) and (A2) and by Lemma 5.2, we have $\lim_{\lambda \to \infty} F_\lambda(z,e) = F_\infty(z,e)$ for all $(z,e)$ in the support of the measure $\hat{N}((0,S] \times dz \times de)$. By the condition (A4), we see that, for any $\lambda \leq \infty$, the integrand $F_\lambda(z,e)$ is dominated by $F_+(z,e)$, which is integrable with respect to the measure $\hat{N}((0,S] \times dz \times de)$. We then appeal to Lebesgue's convergence theorem and obtain $\lim_{\lambda \to \infty} I_\lambda = 0$. Combining this with condition (A3), we obtain the desired result. □

*Remark 6.5.* In the statement of Lemma 6.4, assumption (A4) cannot be removed. For example, let us consider $m_\lambda$ defined by $m_\lambda(x) = x$ for $x \in (0, 1/\lambda)$ and $= x + 1$ for $x \in [1/\lambda, \infty)$, and let $m_\infty(x) = x$. Let $J_\lambda = J_\infty = 0$, $r_\lambda = r_\infty = 0$ and $c_\lambda = c_\infty = c$ for some constant $c > 0$. We then see that all the conditions (A1)–(A3) hold, but we can see (cf. [23]) that $\hat{\eta}_{m_\lambda, J_\lambda, r_\lambda}$ converges in law to $\hat{\eta}_{m_\infty, J_\infty, r_\infty + 1}$, which never coincides in law with $\hat{\eta}_{m_\infty, J_\infty, r_\infty}$.

Let us reduce Propositions 6.1 and 6.2 to Lemma 6.4. For this purpose, we prepare the following lemma.

**Lemma 6.6.** *Let $f$ be a non-negative locally bounded function on $(0, \infty)$. Assume that*

$$f(x) \sim x^\gamma K(x) \qquad \text{as } x \to \infty \tag{6.6}$$

*for some real index $\gamma$ and some slowly varying function $K(x)$ at infinity, and that*

$$\limsup_{x \to 0+} f(x) x^{-\gamma} < \infty. \tag{6.7}$$

*Set $f_\lambda(x) = f(\lambda x)/\{\lambda^\gamma K(\lambda)\}$. Then, for any $\gamma'$ and $\gamma''$ with $\gamma' < \gamma < \gamma''$, there exist constants $C$ and $\lambda_0 > 0$ such that*

$$f_\lambda(x) \leq C \max\{x^{\gamma'}, x^{\gamma''}\} \qquad \text{for all } x > 0 \text{ and all } \lambda > \lambda_0. \tag{6.8}$$

**Proof.** By the assumptions, we may take a constant $C_1$ and a function $\tilde{K}(x)$ defined on $[0,\infty)$ such that the following hold:

(i) $f(x) \leq C_1 x^\gamma \tilde{K}(x)$ for all $x > 0$;
(ii) $\tilde{K}(x)$ is bounded away from 0 and $\infty$ on each compact subset of $[0,\infty)$;
(iii) $\tilde{K}(x)/K(x) \to 1$ as $x \to \infty$ ($\tilde{K}(x)$ is then necessarily slowly varying at $x = \infty$).

We may apply Theorem 1.5.6(ii) of [1], page 25, to the function $\tilde{K}(x)$ and see that there exist constants $C_2$ and $\lambda_0 > 0$ such that

$$\tilde{K}(\lambda x)/\tilde{K}(\lambda) \leq C_2 \max\{x^{\gamma'-\gamma}, x^{\gamma''-\gamma}\} \qquad \text{for all } x > 0 \text{ and all } \lambda > \lambda_0. \tag{6.9}$$

Therefore, we obtain (6.8). □

Thanks to Lemma 6.8, we obtain the following lemma.

**Lemma 6.7.** *Suppose that all the assumptions of either Proposition 6.1 or Proposition 6.2 hold. Condition* (A4) *is then satisfied, where* $J_\infty = V_{(0,d(J))}$ *in the former case and where* $J_\infty = J^{(\beta)}$ *in the latter case.*

**Proof.** Suppose the assumptions of Proposition 6.1 are satisfied. Take numbers $\alpha'$ and $\alpha''$ such that $0 < \alpha' < \alpha < \alpha'' < 1$. Using Lemma 6.6, we know that there exist constants $C$ and $\lambda_0 > 0$ such that $m'_\lambda(x) \leq m'_+(x)$ for all $x > 0$ and all $\lambda > \lambda_0$, where $m_+(x) = C \max\{m^{(\alpha')}(x), m^{(\alpha'')}(x)\}$. Since $J(z) = \infty$ for $z \geq d(J)$, it is obvious that $J_\lambda(z) \geq V_{(0,d(J))}(z)$ for all $z > 0$. Therefore, we may take $J_+ = V_{(0,d(J))}$ to satisfy condition (A4).

Suppose the assumptions of Proposition 6.2 are satisfied. Take numbers $\alpha'$ and $\alpha''$ such that $0 < \alpha' < \alpha < \alpha'' < 1/\beta$. Using Lemma 6.6, we know that there exist constants $C_1$ and $\lambda_1 > 0$ such that $m'_\lambda(x) \leq m'_+(x)$ for all $x > 0$ and all $\lambda > \lambda_1$, where $m_+(x) = C_1 \max\{m^{(\alpha')}(x), m^{(\alpha'')}(x)\}$. Take numbers $\beta'$ and $\beta''$ such that $0 < \beta' < \beta < \beta'' < \min\{1, 1/\alpha''\}$. Using Lemma 6.6 again for $(J_\lambda)^{-1}$, we know that there exist constants $C_2$ and $\lambda_2 > \lambda_1$ such that $J_\lambda(z) \geq J_+(z)$ for all $z > 0$ and all $\lambda > \lambda_2$, where $J_+(z) = C_2 \min\{J^{(\beta')}(z), J^{(\beta'')}(z)\}$. Therefore, we obtain that condition (A4) is satisfied. □

We now proceed to prove Propositions 6.1 and 6.2.

**Proof of Propositions 6.1 and 6.2.** Suppose that all the assumptions of either Proposition 6.1 or Proposition 6.2 hold. By Lemmas 5.3 and 5.5, we know that conditions (A1) and (A2) are satisfied in both cases. It is also obvious that condition (A3) is satisfied for $r_\infty = 0$. By Lemma 6.7, we know that condition (A4) is satisfied. Therefore, the proof follows from Lemma 6.4. □

## 7. Convergence of the Markov process

Propositions 6.1 and 6.2 lead us to the following two propositions, respectively.

**Proposition 7.1 (The convergent case).** *Suppose that the assumptions of Proposition 6.1 are satisfied. The convergence*

$$\hat{X}_{m_\lambda, J_\lambda, r_\lambda} \to \hat{X}_{m^{(\alpha)}, V_{(0,d(J))}, 0} \qquad (J_1) \tag{7.1}$$

*then holds with $\hat{\boldsymbol{P}}$-probability one.*

**Proposition 7.2 (The divergent case).** *Suppose that the assumptions of Proposition 6.2 are satisfied. The convergence*

$$\hat{X}_{m_\lambda, J_\lambda, r_\lambda} \to \hat{X}_{m^{(\alpha)}, J^{(\beta)}, 0} \qquad (J_1) \tag{7.2}$$

*then holds with $\hat{\boldsymbol{P}}$-probability one.*

We introduce the following condition.

(A5) There exist a constant $z_0 > 0$ and a right-continuous non-decreasing function $J_+ : (0, \infty) \to [0, \infty]$ with

$$\int_{z_0}^\infty \frac{\mathrm{d}z}{J_+(z)} < \infty \tag{7.3}$$

such that $J_\lambda(z) \geq J_+(z)$ for all $z > z_0$.

We now obtain the following continuity lemma for the Markov process.

**Lemma 7.3.** *Suppose that conditions* (A1)–(A3) *and* (A5) *hold and that the convergence*

$$\lim_{\lambda \to \infty} \sup_{s \in [0,S]} |\hat{\eta}_{m_\lambda, J_\lambda, r_\lambda}(s) - \hat{\eta}_{m_\infty, J_\infty, r_\infty}(s)| = 0 \qquad \text{for all } S > 0 \tag{7.4}$$

*holds with $\hat{\boldsymbol{P}}$-probability one. The convergence*

$$\hat{X}_{m_\lambda, J_\lambda, r_\lambda} \to \hat{X}_{m_\infty, J_\infty, r_\infty} \qquad (J_1) \tag{7.5}$$

*then holds with $\hat{\boldsymbol{P}}$-probability one.*

**Proof.** 1. Since $\boldsymbol{n}_{\mathrm{BE}}(M > x) = 1/x$, we have

$$\int_{(0,\infty) \times E} F(z, e) 1_{\{M(e) > \varepsilon\}} \, \mathrm{d}z \otimes \boldsymbol{n}_{\mathrm{BE}}(\mathrm{d}e) < \infty \qquad \text{for all } \varepsilon > 0, \tag{7.6}$$

where $F(z, e) = 1_{\{0 < z \leq z_0\}} + 1_{\{z > z_0, M(e) > J_+(z)\}}$. In fact, the left-hand side of (7.6) is dominated by

$$\frac{z_0}{\varepsilon} + \int_{z_0}^\infty \frac{1}{\max\{J_+(z), \varepsilon\}} \, \mathrm{d}z, \tag{7.7}$$

which turns out to be finite by the assumption (7.3) of (A5). Hence, we obtain that

$$\int_{(0,\infty)\times E} F(z,e)1_{\{M(e)>\varepsilon\}}\hat{\boldsymbol{N}}((0,s]\times \mathrm{d}z\times \mathrm{d}e) < \infty \qquad \text{for all } s\geq 0 \text{ and } \varepsilon>0 \quad (7.8)$$

holds with $\hat{\boldsymbol{P}}$-probability one. In addition, recall that we can apply Lemma 5.2 in this case and obtain that

$$\lim_{\lambda\to\infty}\sup_{t\geq 0}|e_{m_\lambda,J_\lambda(z)}(t) - e_{m_\infty,J_\infty(z)}(t)| = 0 \tag{7.9}$$

for all $(z,e)$ in the support of $\hat{\boldsymbol{N}}((0,\infty)\times \mathrm{d}z\times \mathrm{d}e)$

holds with $\hat{\boldsymbol{P}}$-probability one. Thus, there exists $\hat{\Omega}^* \in \hat{\mathcal{F}}$ with $\hat{\boldsymbol{P}}(\hat{\Omega}^*)=1$ on which (7.8), (7.9) and (7.4) hold. Let $\hat{\omega}\in\hat{\Omega}^*$ be fixed until the end of the proof.

2. We shall construct a family of functions $\{\Lambda_\lambda : \lambda > 0\}$ (which may depend on $\hat{\omega}$) imitating Stone [23]. For any $\varepsilon > 0$, the support of the point process

$$F(z,e)1_{\{M>\varepsilon\}}(e)\hat{\boldsymbol{N}}(\mathrm{d}s\times \mathrm{d}z\times \mathrm{d}e) \tag{7.10}$$

on $(0,\infty)\times(0,\infty)\times E^*$ is enumerated by $\{(s^{\varepsilon,(i)}, z^{\varepsilon,(i)}, e^{\varepsilon,(i)}): i=1,2,\ldots\}$ such that $s^{\varepsilon,(1)} < s^{\varepsilon,(2)} < \cdots$. Define

$$\Lambda_{\varepsilon,\lambda}(\hat{\eta}_{m_\infty,J_\infty,r_\infty}(s^{\varepsilon,(i)}-)) = \hat{\eta}_{m_\lambda,J_\lambda,r_\lambda}(s^{\varepsilon,(i)}-), \qquad i=1,2,\ldots, \tag{7.11}$$

$$\Lambda_{\varepsilon,\lambda}(\hat{\eta}_{m_\infty,J_\infty,r_\infty}(s^{\varepsilon,(i)})) = \hat{\eta}_{m_\lambda,J_\lambda,r_\lambda}(s^{\varepsilon,(i)}), \qquad i=1,2,\ldots, \tag{7.12}$$

and extend $\Lambda_{\varepsilon,\lambda}$ to a continuous function on $(0,\infty)$ by linear interpolation. If the number $n$ of $s^{\varepsilon,(i)}$'s is finite, then we set $\Lambda_{\varepsilon,\lambda}(t) = t - t_n + \Lambda_{\varepsilon,\lambda}(t_n)$ for $t > t_n := \hat{\eta}_{m_\infty,J_\infty,r_\infty}(s^{\varepsilon,(n)})$. We define $\Lambda_\lambda = \Lambda_{1/\lambda,\lambda}$. Since $\hat{\eta}_{m_\lambda,J_\lambda,r_\lambda}(\infty) = \hat{\eta}_{m_\infty,J_\infty,r_\infty}(\infty) = \infty$, we see that $\Lambda_\lambda(\infty) = \infty$ and hence $\Lambda_\lambda$ is a homeomorphism of $[0,\infty)$. Since (7.4) holds, it is immediate that

$$\lim_{\lambda\to\infty}\sup_{t\in[0,T]}|\Lambda_\lambda(t) - t| = 0 \tag{7.13}$$

for all $T > 0$.

3. Let $\varepsilon > 0$ be fixed. It suffices to show that

$$\limsup_{\lambda\to\infty}\sup_{t\in[0,T]}|\hat{X}_{m_\lambda,J_\lambda,r_\lambda}(\Lambda_\lambda(t)) - \hat{X}_{m_\infty,J_\infty,r_\infty}(t)| \leq 2\varepsilon \tag{7.14}$$

for all $T > 0$.

For $\lambda$ such that $1/\varepsilon < \lambda \leq \infty$, we set

$$I_\lambda^{\varepsilon,(i)} = [\hat{\eta}_{m_\lambda,J_\lambda,r_\lambda}(s^{\varepsilon,(i)}-), \hat{\eta}_{m_\lambda,J_\lambda,r_\lambda}(s^{\varepsilon,(i)})) \subset (0,\infty), \qquad i=1,2,\ldots. \tag{7.15}$$

By definition, we have $\Lambda_\lambda(I_\infty^{\varepsilon,(i)}) = I_\lambda^{\varepsilon,(i)}$.

4. Let $t \notin \bigcup_i I_\infty^{\varepsilon,(i)}$. We then have $\Lambda_\lambda(t) \notin \bigcup_i I_\lambda^{\varepsilon,(i)}$ for all $\lambda > 1/\varepsilon$. For $1/\varepsilon < \lambda \leq \infty$, we take $(s_\lambda, z_\lambda, e_\lambda)$ such that $\hat{\eta}_{m_\lambda, J_\lambda, r_\lambda}(s_\lambda-) \leq \Lambda_\lambda(t) < \hat{\eta}_{m_\lambda, J_\lambda, r_\lambda}(s_\lambda)$, if it exists, where $\Lambda_\infty(t) = t$. If such a point $(s_\lambda, z_\lambda, e_\lambda)$ does not exist, then $\hat{X}_{m_\lambda, J_\lambda, r_\lambda}(t) = 0$. If $(s_\lambda, z_\lambda, e_\lambda)$ exists, then we have $M(e_\lambda) \leq \varepsilon$. In fact, if, in addition, $z_\lambda > z_0$, then $M(e_\lambda) \geq J_\lambda(z_\lambda) \geq J_+(z_\lambda)$, by assumption (A5). In both cases, we have $\hat{X}_{m_\lambda, J_\lambda, r_\lambda}(\Lambda_\lambda(t)) \leq \varepsilon$. Hence, we obtain

$$\sup_{t \notin \bigcup_i I^{\varepsilon,(i)}} |\hat{X}_{m_\lambda, J_\lambda, r_\lambda}(\Lambda_\lambda(t)) - \hat{X}_{m_\infty, J_\infty, r_\infty}(t)| < 2\varepsilon \qquad \text{for all } \lambda > 1/\varepsilon. \tag{7.16}$$

Let $t \in I_\infty^{\varepsilon,(i)}$ for some $i$. Write $(s^{\varepsilon,(i)}, z^{\varepsilon,(i)}, e^{\varepsilon,(i)})$ simply as $(s, z, e)$ for now. We then have

$$\hat{X}_{m_\infty, J_\infty, r_\infty}(t) = e_{m_\infty, J_\infty(z)}(t - \hat{\eta}_{m_\infty, J_\infty, r_\infty}(s-)) \tag{7.17}$$

and, since $\Lambda_\lambda(t) \in I_\lambda^{\varepsilon,(i)}$, we have

$$\hat{X}_{m_\lambda, J_\lambda, r_\lambda}(\Lambda_\lambda(t)) = e_{m_\lambda, J_\lambda(z)}(\Lambda_\lambda(t) - \hat{\eta}_{m_\lambda, J_\lambda, r_\lambda}(s-)). \tag{7.18}$$

Since (7.9), (7.4) and (7.13) hold, we obtain

$$\lim_{\lambda \to \infty} \hat{X}_{m_\lambda, J_\lambda, r_\lambda}(\Lambda_\lambda(t)) = \hat{X}_{m_\infty, J_\infty, r_\infty}(t). \tag{7.19}$$

5. Since we have at most a finite number of $i$'s such that $I_\infty^{\varepsilon,(i)} \cap [0,T] \neq \varnothing$, it follows from (7.16) and (7.19) that (7.14) holds for all $T > 0$. We now conclude that (6.4) holds. □

We now prove Propositions 7.1 and 7.2.

**Proof of Propositions 7.1 and 7.2.** Suppose that all the assumptions of either Propositions 6.1 or Proposition 6.2 hold. We then know that all the conditions (A1)–(A4) hold. It is obvious that condition (A4) implies condition (A5). Therefore, the desired result follows from Lemma 7.3. □

## 8. Removal of the extra assumptions

We remove the extra assumptions (6.1) from Proposition 6.1 and (6.3) from Proposition 6.2 and obtain in-probability continuity results as follows.

**Proposition 8.1 (The convergent case).** *Suppose that the assumptions of Theorem 2.5 hold. Set $v(\lambda) = \lambda$ and adopt the notation (4.3) and (4.4). Then,*

$$\lim_{\lambda \to \infty} \sup_{s \in [0,S]} |\hat{\eta}_{m_\lambda, J_\lambda, r_\lambda}(s) - \hat{\eta}_{m^{(\alpha)}, V_{(0, d(J))}, 0}(s)| = 0 \tag{8.1}$$

*and*

$$\hat{X}_{m_\lambda, J_\lambda, r_\lambda} \to \hat{X}_{m^{(\alpha)}, V_{(0,d(J))}, 0} \qquad (J_1) \tag{8.2}$$

*hold in probability for all $S > 0$.*

**Proposition 8.2 (The divergent case).** *Suppose that the assumptions of Theorem 2.6 hold. Set $v(\lambda) = \lambda^\beta / L(\lambda)$ and adopt the notation (4.3) and (4.4). Then,*

$$\lim_{\lambda \to \infty} \sup_{s \in [0,S]} |\hat{\eta}_{m_\lambda, J_\lambda, r_\lambda}(s) - \hat{\eta}_{m^{(\alpha)}, J^{(\beta)}, 0}(s)| = 0 \tag{8.3}$$

*and*

$$\hat{X}_{m_\lambda, J_\lambda, r_\lambda} \to \hat{X}_{m^{(\alpha)}, J^{(\beta)}, 0} \qquad (J_1) \tag{8.4}$$

*hold in probability for all $S > 0$.*

Theorems 2.5 and 2.6 immediately follow from Propositions 8.1 and 8.2, respectively.

In order to prove Propositions 8.1 and 8.2, we prepare two lemmas. The first one is the following.

**Lemma 8.3.** *Let $(m, J)$ be a pair which satisfies condition* (C). *Assume that $\mathrm{d}m(x) = 0$ on $(x_0, \infty)$ for some $x_0 > 0$. Then, for any $\gamma < 1$,*

$$\lim_{\lambda \to \infty} \frac{1}{\lambda^{1/\gamma}} \hat{\eta}_{m, J, 0}(\lambda s) = 0 \tag{8.5}$$

*holds in probability.*

**Proof.** Taking a Laplace transform, we can see that it suffices to show that

$$\lim_{\varepsilon \to 0+} \varepsilon^{-\gamma} \left\{ c(J) \boldsymbol{n}_{\mathrm{BE}}[1 - \mathrm{e}^{-\varepsilon \zeta(e_m)}] + \int_{(0,\infty)} \frac{\mathrm{d}z}{J(z)} \boldsymbol{Q}_{\mathrm{BM}}^{J(z)}[1 - \mathrm{e}^{-\varepsilon \zeta(e_m)}] \right\} = 0. \tag{8.6}$$

It is well known that

$$\boldsymbol{Q}_{\mathrm{BM}}^{x}[1 - \mathrm{e}^{-\varepsilon \zeta(e_m)}] = 1 - g_\varepsilon(x), \tag{8.7}$$

where $g_\varepsilon(x)$ satisfies

$$1 - g_\varepsilon(x) = \varepsilon \int_0^x \mathrm{d}y \int_{(y, x_0]} g_\varepsilon(z) \, \mathrm{d}m(z) \tag{8.8}$$

and $g_\varepsilon(x) = g_\varepsilon(x_0)$ for all $x > x_0$. We use the inequality $g_\varepsilon \leq 1$ to obtain

$$1 - g_\varepsilon(x) \leq \varepsilon \int_0^{\min\{x, x_0\}} m((y, x_0]) \, \mathrm{d}y \qquad \text{for all } x > 0. \tag{8.9}$$

Hence, we obtain

$$\varepsilon^{-1} \int_{(0,\infty)} \frac{\mathrm{d}z}{J(z)} \boldsymbol{Q}_{\mathrm{BM}}^{J(z)}[1 - \mathrm{e}^{-\varepsilon\zeta(e_m)}] \leq \int_{(0,\infty)} \frac{\mathrm{d}z}{J(z)} \int_0^{\min\{J(z),x_0\}} m((y,x_0])\,\mathrm{d}y. \quad (8.10)$$

The right-hand side turns out to be finite by assumption (3.16).

Suppose that $c(J) > 0$. The origin for $\mathcal{L}_m$ must then be exit and entrance, that is, $m(0+)$ is finite. Since $\boldsymbol{n}_{\mathrm{BE}}[1 - \mathrm{e}^{-\varepsilon\zeta(e_m)}] = \lim_{x \to 0+} \frac{1}{x}\boldsymbol{Q}_{\mathrm{BM}}^x[1 - \mathrm{e}^{-\varepsilon\zeta(e_m)}]$, we know that

$$\varepsilon^{-1}\boldsymbol{n}_{\mathrm{BE}}[1 - \mathrm{e}^{-\varepsilon\zeta(e_m)}] = \int_{(0,x_0]} g_\varepsilon(z)\,\mathrm{d}m(z) \leq m((0,x_0]) < \infty. \quad (8.11)$$

Therefore, the proof is complete. $\square$

**Lemma 8.4.** *Let $(m, J)$ be a pair which satisfies condition (C). Suppose that $m(x)$ satisfies (M) for $\alpha \in (0,\infty)$ and $x_0 > 0$ and that $\mathrm{d}m(x) = 0$ on $(0, x_0)$. Suppose that $d(J) < \infty$ and that $c(J) = 0$ when $\alpha \geq 1$. Then, for any $\gamma < \min\{1, \alpha\}$,*

$$\lim_{\lambda \to \infty} \frac{1}{\lambda^{1/\gamma}} \hat{\eta}_{m,J,0}(\lambda s) = 0 \quad (8.12)$$

*holds in probability.*

**Proof.** 1. Consider the case where $\alpha < 1$. For any $\nu \in (\gamma, \alpha)$, there exists a constant $C_1$ such that $m'(x) \leq C_1 m^{(\nu)}(x)$ for all $x > 0$. Since

$$\hat{\eta}_{m,J,0}(\lambda s) \leq C_1 \hat{\eta}_{m^{(\nu)},J,0}(\lambda s), \quad (8.13)$$

it suffices to show that

$$\lim_{\varepsilon \to 0+} \varepsilon^{-\gamma}\left\{c(J)\boldsymbol{n}_{\mathrm{BE}}[1 - \mathrm{e}^{-\varepsilon\zeta(e_{m^{(\nu)}})}] + \int_{(0,d(J))} \frac{\mathrm{d}z}{J(z)} \boldsymbol{Q}_{\mathrm{BM}}^{J(z)}[1 - \mathrm{e}^{-\varepsilon\zeta(e_{m^{(\nu)}})}]\right\} = 0. \quad (8.14)$$

It is well known that

$$\boldsymbol{Q}_{\mathrm{BM}}^x[1 - \mathrm{e}^{-\varepsilon\zeta(e_{m^{(\nu)}})}] = 1 - g(\varepsilon^\nu x), \quad (8.15)$$

where $g(x)$ satisfies

$$1 - g(x) = -g'(x_0)x + \int_0^x \mathrm{d}y \int_y^{x_0} g(z)\,\mathrm{d}m^{(\nu)}(z), \quad (8.16)$$

and that

$$\boldsymbol{n}_{\mathrm{BE}}[1 - \mathrm{e}^{-\varepsilon\zeta(e_{m^{(\nu)}})}] = \varepsilon^\nu\left\{-g'(x_0) + \int_0^{x_0} g(z)\,\mathrm{d}m^{(\nu)}(z)\right\}. \quad (8.17)$$

Since $m^{(\nu)}((0, x_0)) < \infty$ and $g(z) \leq 1$ for all $z > 0$, there exists a constant $C$ such that $\varepsilon^{-\nu} \boldsymbol{Q}_{\mathrm{BM}}^x [1 - \mathrm{e}^{-\varepsilon \zeta(e_{m^{(\nu)}})}] \leq C$ and $\varepsilon^{-\nu} \boldsymbol{n}_{\mathrm{BE}}[1 - \mathrm{e}^{-\varepsilon \zeta(e_{m^{(\nu)}})}] \leq C$ for all $\varepsilon > 0$. Therefore, we obtain (8.14).

2. In the case where $\alpha = 1$, we can prove the desired convergence in almost the same way as 1. The only difference is to use $c(J) = 0$. We omit the details.

3. Consider the case where $\alpha > 1$. Then, $c(J) = 0$. Taking a Laplace transform, it suffices to show that

$$\lim_{\varepsilon \to 0+} \varepsilon^{-\gamma} \int_{(0,d(J))} \frac{\mathrm{d}z}{J(z)} \boldsymbol{Q}_{\mathrm{BM}}^{J(z)} [1 - \mathrm{e}^{-\varepsilon \zeta(e_m)}] = 0. \quad (8.18)$$

It is well known that $\boldsymbol{Q}_{\mathrm{BM}}^x [1 - \mathrm{e}^{-\varepsilon \zeta(e_m)}] = 1 - g_\varepsilon(x)$, where $g_\varepsilon(x)$ satisfies

$$1 - g_\varepsilon(x) = \varepsilon \int_0^x \mathrm{d}y \int_y^\infty g_\varepsilon(z) \, \mathrm{d}m(z). \quad (8.19)$$

We use the inequality $g_\varepsilon \leq 1$ to obtain

$$\varepsilon^{-1} \int_{(0,d(J))} \frac{\mathrm{d}z}{J(z)} \boldsymbol{Q}_{\mathrm{BM}}^{J(z)} [1 - \mathrm{e}^{-\varepsilon \zeta(e_m)}] \leq \int_{(0,d(J))} \frac{\mathrm{d}z}{J(z)} \int_0^{J(z)} m((y, \infty)) \, \mathrm{d}y. \quad (8.20)$$

The right-hand side is finite by condition (C) and we therefore obtain (8.18). □

**Proof of Propositions 8.1 and 8.2.** 1. In the case of Proposition 8.1, we take $z_0 = d(J)$, $J_+ = V_{(0,d(J))}$ and $(m_\infty, J_\infty, r_\infty) = (m^{(\alpha)}, V_{(0,d(J))}, 0)$. In the case of Proposition 8.2, using Lemma 6.6, we have that there exist constants $C > 0$, $z_0 > 0$ and $\beta'$ with $\beta < \beta' < \max\{1, 1/\alpha\}$ such that $J_\lambda(z) \geq J_+(z)$ for all $z > z_0$ and $\lambda > 0$, where $J_+(z) = CJ^{(\beta')}(z)$. In this case, we take $(m_\infty, J_\infty, r_\infty) = (m^{(\alpha)}, J^{(\beta)}, 0)$. Now, in both cases, the triplet $(m_\lambda^0, J_\lambda^0, r_\lambda)$ satisfies all of the assumptions (A1)–(A3) and (A5).

Let us define $m^0$ and $m^1$ by

$$m^0(x) = m(\max\{x, x_0\}) \quad \text{and} \quad m^1(x) = m(\min\{x, x_0\}). \quad (8.21)$$

Let us define $J^0$ and $J^1$ by

$$J^0(z) = J(z + z_0) \quad \text{and} \quad J^1(z) = \begin{cases} J(z), & \text{on } (0, z_0), \\ \infty, & \text{on } [z_0, \infty). \end{cases} \quad (8.22)$$

We define $m_\lambda^0$ and $m_\lambda^1$ (resp. $J_\lambda^0$ and $J_\lambda^1$) in the same way as $m_\lambda$ in (4.1) (resp. $J_\lambda$ in (4.4)). The triplet $(m_\lambda^0, J_\lambda^0, r_\lambda)$ then satisfies all of the assumptions (A1)–(A4). We now have

$$\hat{\eta}_{m_\lambda, J_\lambda, r_\lambda} = \hat{\eta}_\lambda^0 + \hat{\eta}_\lambda^1 + \hat{\eta}_\lambda^2, \quad (8.23)$$

where

$$\hat{\eta}_\lambda^1 = \hat{\eta}_{m_\lambda^0, J_\lambda^1, 0}, \qquad \hat{\eta}_\lambda^2 = \hat{\eta}_{m_\lambda^1, J_\lambda, 0} \quad (8.24)$$

and

$$\hat{\eta}^0_\lambda(s) = r_\lambda + \int_{[v(\lambda)/\lambda z_0, d(J)) \times E} \zeta(e_{m_\lambda, J_\lambda(z)}) 1_{\{M(e) \geq J_\lambda(z)\}} \hat{\boldsymbol{N}}((0,s] \times \mathrm{d}z \times \mathrm{d}e) \quad (8.25)$$

for $s \geq 0$. By the translation invariance in $z$ of the characteristic measure of $\hat{\boldsymbol{N}}(\mathrm{d}s \times \mathrm{d}z \times \mathrm{d}e)$, it is obvious that

$$\hat{\eta}^0_\lambda \stackrel{\text{law}}{=} \hat{\eta}_{m^0_\lambda, J^0_\lambda, r_\lambda}. \quad (8.26)$$

Since the triplet $(m^0, J^0, 0)$ satisfies the assumptions, we may apply Lemma 6.4 and obtain

$$\lim_{\lambda \to \infty} \sup_{s \in [0,S]} |\hat{\eta}^0_\lambda(s) - \hat{\eta}_{m_\infty, J_\infty, r_\infty}(s)| = 0 \qquad \text{for all } S > 0 \quad (8.27)$$

$\hat{\boldsymbol{P}}$-almost surely.

Using Lemma 4.1 again, we have

$$\hat{\eta}^1_\lambda(s) \stackrel{\text{law}}{=} \frac{1}{u(\lambda)} \hat{\eta}_{m^0, J^1, 0}(v(\lambda)s) \quad \text{and} \quad \hat{\eta}^2_\lambda(s) \stackrel{\text{law}}{=} \frac{1}{u(\lambda)} \hat{\eta}_{m^1, J, 0}(v(\lambda)s). \quad (8.28)$$

We note that the pair $(m^0, J^1)$ (resp. $(m^1, J)$) satisfies the assumptions of Lemma 8.4 (resp. Lemma 8.3). In the case of Proposition 8.1, we have $u(v^{-1}(\lambda)) \sim \lambda^{1/\alpha} K(\lambda)$ as $\lambda \to \infty$. In the case of Proposition 8.2, we have $u(v^{-1}(\lambda)) \sim \lambda^{1/(\alpha\beta)} \tilde{K}(\lambda)$ as $\lambda \to \infty$ for some slowly varying function $\tilde{K}$ at infinity, where $v^{-1}$ is an asymptotic inverse of $v$. In both cases, we have $u(v^{-1}(\lambda)) \sim \lambda^{1/\gamma}$ for some $\gamma < 1$. Hence, by Lemmas 8.4 and 8.3, we obtain

$$\lim_{\lambda \to \infty} \hat{\eta}^1_\lambda(s) = 0 \quad \text{and} \quad \lim_{\lambda \to \infty} \hat{\eta}^2_\lambda(s) = 0 \quad (8.29)$$

in probability, for all $s > 0$.

Consequently, we obtain

$$\lim_{\lambda \to \infty} \sup_{s \in [0,S]} |\hat{\eta}_{m_\lambda, J_\lambda, r_\lambda}(s) - \hat{\eta}_{m_\infty, J_\infty, r_\infty}(s)| = 0, \quad (8.30)$$

in probability, for all $S > 0$. Let $\lambda(n)$ be an arbitrary sequence of $(0, \infty)$ such that $\lambda(n) \to \infty$. We can then take a subsequence $\lambda(n_k)$ along which (8.30) holds for $S > 0$ with $\hat{\boldsymbol{P}}$-probability one. We may now apply Lemma 7.3 to obtain

$$X_{m_\lambda, J_\lambda, r_\lambda} \to X_{m_\infty, J_\infty, r_\infty} \qquad (J_1) \quad (8.31)$$

along the subsequence $\lambda = \lambda(n_k)$. This means that the convergence (8.31) occurs in probability. Therefore, we obtain the desired conclusions. $\square$

# Acknowledgements

The author was supported by JSPS Research Fellowships for Young Scientists. He would like to thank Professor Masatoshi Fukushima for granting him access to Professor Kiyosi Itô's unpublished lecture note [8].